\documentclass[12pt]{amsart}

\usepackage{amssymb}
\usepackage[colorlinks]{hyperref}
\usepackage[margin=1.25in]{geometry}
\usepackage{enumitem}
\usepackage{tikz}
\usetikzlibrary{matrix,arrows}
\usepackage{verbatim}

\theoremstyle{definition}

\newtheorem*{definition}{Definition}

\author{Federico Galetto}
\title{Setting the scene for Betti characters}
\date{\today}

\address{Department of Mathematics and Statistics, Cleveland State University, Cleveland, OH, 44115-2215, USA}
\email{\href{mailto:f.galetto@csuohio.edu}{\nolinkurl{f.galetto@csuohio.edu}}}
\urladdr{\href{https://math.galetto.org}{\nolinkurl{https://math.galetto.org}}}



\begin{document}

\begin{abstract}
  Finite group actions on free resolutions and modules arise naturally
  in many interesting examples.  Understanding these actions amounts
  to describing the terms of a free resolution or the graded
  components of a module as group representations which, in the non
  modular case, are completely determined by their characters.  With
  this goal in mind, we introduce a Macaulay2 package for computing
  characters of finite groups on free resolutions and graded
  components of finitely generated graded modules over polynomial
  rings.
\end{abstract}

\maketitle{}

\section{Introduction}
\label{sec:introduction}

Let $R$ be a polynomial ring over a field $\Bbbk$, and $M$ a finitely
generated graded $R$-module. Let $G$ be a linearly reductive group
acting $\Bbbk$-linearly on $R$ and $M$, and assume these actions
preserve degrees and distribute over $R$-multiplication. If
$F_\bullet$ is a minimal graded free resolution of $M$, then the
action of $G$ extends to $F_\bullet$.  More precisely, if
$\mathfrak{m}$ denotes the irrelevant maximal ideal of $R$, then the
finite dimensional vector spaces $F_i / \mathfrak{m} F_i$ carry a
natural structure of graded $G$-representations (see \cite[Proposition
2.4.9, Remark 2.4.10]{MR3424030} for details).  This additional
structure makes a resolution more rigid as the differentials must
commute with the group action; in some cases, this makes it possible
to construct the differentials explicitly using representation theory
(see, for example, \cite{MR2878669,MR2769324}). Understanding how $G$
acts on the modules $F_i / \mathfrak{m} F_i$ may also lead to
interesting combinatorial descriptions of the Betti numbers of $M$,
such as in \cite[Corollary 4.12]{MR4105544}.  Free resolutions
equipped with group actions (also known as equivariant resolutions)
have found many important applications, such as the computation of
Betti numbers of determinantal varieties \cite{MR520233}, and a proof
of the existence of pure free resolutions \cite{MR2918721} (a central
aspect of Boij-S\"oderberg theory).

From a computational perspective, the Macaulay2 \cite{M2} package
\texttt{HighestWeights} \cite{MR3438710} allows users to determine the
representation theoretic structure of an equivariant resolution with
the action of a semisimple Lie group in characteristic zero. Recent
publications
\cite{MR3215587,MR3779986,MR3799202,MR4043827,MR4105544,MR4108337,MR4155201,2010.06522,2012.13732,MR4136641}
point to an interest in equivariant resolutions with actions of finite
groups, particularly symmetric groups. However, at the time of
writing, no software solution is available to compute such
actions. The present article introduces the Macaulay2 package
\texttt{BettiCharacters}\footnote{Available at
  \url{https://github.com/galettof/BettiCharacters}.} to fill this
gap.  In the non modular case (i.e., when the characteristic of the
field does not divide the order of the group), finite dimensional
representations of finite groups are determined, up to isomorphism, by
their characters (see \cite[Chapter 2]{MR0450380} for an introduction
to the subject).  Thus understanding the $G$-action on a minimal free
resolution $F_\bullet$ amounts to describing the graded characters of
the representations $F_i / \mathfrak{m} F_i$ or, equivalently, the
characters of the graded components $(F_i / \mathfrak{m} F_i)_j$.  The
uniqueness of minimal free resolutions implies
$(F_i / \mathfrak{m} F_i)_j \cong \operatorname{Tor}^R_i (M, \Bbbk)_j$
as $G$-representations.  Moreover, the character of $G$ on
$\operatorname{Tor}^R_i (M, \Bbbk)_j$ evaluated at the identity of $G$
is the dimension of $\operatorname{Tor}^R_i (M, \Bbbk)_j$ as a
$\Bbbk$-vector space, i.e., the $(i,j)$-th Betti number of $M$.
Therefore we adopt the following definition, after which the package
is named.
\begin{definition}
  The \emph{$(i,j)$-th Betti character of $G$ on $M$}, denoted
  $\beta^G_{i,j} (M)$, is the character of $G$ on
  $\operatorname{Tor}^R_i (M, \Bbbk)_j$.
\end{definition}
The package \texttt{BettiCharacters} implements the algorithm
described in \cite[Algorithm 1]{2106.14071}, which in essence
propagates the trace of the group action from the module $M$ through a
(previously computed) minimal free resolution $F_\bullet$. In
addition, \texttt{BettiCharacters} also allows users to compute the
characters of $G$ on the graded components of $M$.  The rest of this
article illustrates the main functionalities of the
package\footnote{All computations performed in Macaulay2 1.18 on
  Debian 10, using \texttt{BettiCharacters} 1.0.}.

\section{Example: a symmetric shifted ideal}
\label{sec:example-1:-symmetric}

Consider the ideal generated by all quadratic squarefree monomials in
a ring with four variables.

\begin{verbatim}
i1 : R = QQ[x_1..x_4];

i2 : I = ideal apply(subsets(gens R,2),product)

o2 = ideal (x x , x x , x x , x x , x x , x x )
             1 2   1 3   2 3   1 4   2 4   3 4

o2 : Ideal of R

i3 : RI = res I

      1      6      8      3
o3 = R  <-- R  <-- R  <-- R  <-- 0
                                  
     0      1      2      3      4

o3 : ChainComplex
\end{verbatim}

The symmetric group on four elements acts by permuting the ring
variables and, in doing so, preserves the ideal. Thus the action
passes to the quotient and its minimal free resolution.  The
equivariant structure of the resolution is described in \cite[Theorem
4.11]{MR4105544} and \cite[Theorem 4.1]{MR3779986}.  The ideal also
belongs to the larger class of symmetric shifted ideals, whose
equivariant resolutions are described in \cite[Theorem
6.2]{MR4108337}.  To verify these results computationally, we first
define the group action. Since characters are class functions (i.e.,
constant on conjugacy classes) it is enough to define a single group
element per conjugacy class. In this case, the elements can be
conveniently introduced as one-row matrices of substitutions for the
ring variables. Then we set up the action on the resolution as an
object of type \texttt{ActionOnComplex} using the \texttt{action}
method.

\begin{verbatim}
i4 : G = {matrix{{x_2,x_3,x_4,x_1}}, matrix{{x_2,x_3,x_1,x_4}},
         matrix{{x_2,x_1,x_4,x_3}}, matrix{{x_2,x_1,x_3,x_4}},
         matrix{{x_1,x_2,x_3,x_4}} }

o4 = {| x_2 x_3 x_4 x_1 |, | x_2 x_3 x_1 x_4 |, | x_2 x_1 x_4 x_3 |,
     ---------------------------------------------------------------
     | x_2 x_1 x_3 x_4 |, | x_1 x_2 x_3 x_4 |}

o4 : List

i5 : needsPackage "BettiCharacters";

i6 : A = action(RI,G)

o6 = ChainComplex with 5 actors

o6 : ActionOnComplex
\end{verbatim}

Now we can use the \texttt{character} method to compute Betti
characters.

\begin{verbatim}
i7 : character A

o7 = HashTable{0 => GradedCharacter{{0} => Character{1, 1, 1, 1, 1}}  }
               1 => GradedCharacter{{2} => Character{0, 0, 2, 2, 6}}
               2 => GradedCharacter{{3} => Character{0, -1, 0, 0, 8}}
               3 => GradedCharacter{{4} => Character{1, 0, -1, -1, 3}}

o7 : HashTable
\end{verbatim}

The output is a \texttt{HashTable} with integers (corresponding to
homological degrees) as keys and objects of type
\texttt{GradedCharacter} as values. A \texttt{GradedCharacter} is also
a \texttt{HashTable} with (multi)degrees as keys and objects of type
\texttt{Character} as values. Finally, a \texttt{Character} is simply
a list of traces of the previously defined group elements, and is
meant to be read as a row of a character table.  The methods in
\texttt{BettiCharacters} are completely independent of the
group. However, we invite the reader to consult Examples 1 and 2 in
the documentation to see how this package can be used in conjunction
with \texttt{SpechtModule} \cite{SpechtModuleSource} to decompose
characters of a symmetric group into irreducible characters.

\section{Example: Klein point configuration}
\label{sec:example-2:-klein}

We consider the Klein configuration of points in the projective
plane. The defining ideal $I$ is explicitly constructed in
\cite[Proposition 7.3]{MR4043827}. Although $I$ is defined over the
rationals, we work over the cyclomotic field obtained by adjoining a
primitive seventh root of unity for the purpose of defining a group
action.

\begin{verbatim}
i1 : kk=toField(QQ[a]/ideal(sum apply(7,i->a^i)));

i2 : R=kk[x,y,z];

i3 : f4=x^3*y+y^3*z+z^3*x

      3     3       3
o3 = x y + y z + x*z

o3 : R

i4 : f6=-1/54*det(jacobian transpose jacobian f4)

        5    5      2 2 2      5
o4 = x*y  + x z - 5x y z  + y*z

o4 : R

i5 : I=minors(2,jacobian matrix{{f4,f6}});

o5 : Ideal of R
\end{verbatim}

The unique simple group $G$ of order 168 acts on the projective plane
preserving the Klein configuration. This induces an action on our
polynomial ring preserving the ideal $I$. The action (which is
minimally defined over our cyclomotic field) is explicitly described
in \cite[\S 2.2]{MR4043827}. In particular, the group is generated by
elements \texttt{g} of order 7, \texttt{h} of order 3, and \texttt{i}
of order 2. Since we are interested in some characters of $G$, we need
a representative for each conjugacy class; therefore, in addition to
\texttt{g}, \texttt{h}, and \texttt{i}, we also consider the identity
element, the inverse of \texttt{g}, and an element \texttt{j} of order
4. We define all these group elements as matrices.

\begin{verbatim}
i6 : g=matrix{{a^4,0,0},{0,a^2,0},{0,0,a}};

              3        3
o6 : Matrix kk  <--- kk

i7 : h=matrix{{0,1,0},{0,0,1},{1,0,0}};

              3        3
o7 : Matrix ZZ  <--- ZZ

i8 : i=(2*a^4+2*a^2+2*a+1)/7 * matrix{{a-a^6,a^2-a^5,a^4-a^3},
         {a^2-a^5,a^4-a^3,a-a^6}, {a^4-a^3,a-a^6,a^2-a^5}};

              3        3
o8 : Matrix kk  <--- kk

i9 : j=-1/(2*a^4+2*a^2+2*a+1) * matrix{{a^5-a^4,1-a^5,1-a^3},
         {1-a^5,a^6-a^2,1-a^6}, {1-a^3,1-a^6,a^3-a}};

              3        3
o9 : Matrix kk  <--- kk

i10 : G={id_(R^3),i,h,j,g,inverse g};
\end{verbatim}

As proved in \cite[Theorem 4.4]{MR3351566} and \cite[Proposition
8.1]{MR4043827}, the symbolic cube $I^{(3)}$ is not contained in the
square $I^2$. The second proof of \cite[Proposition 8.1]{MR4043827}
reduces the failure of containment to showing the graded component of
degree 21 in the quotient $I^{(2)} / I^2$ is a trivial
$G$-representation. By local duality, this is equivalent to showing
that the last module in a minimal free resolution of $I^2$ is
generated in degree 24 by a one dimensional trivial $G$-module.  We
proceed to compute the character of $G$ on the last module of the
resolution of $I^2$.

\begin{verbatim}
i11 : I2=I^2;

o11 : Ideal of R

i12 : RI2=res I2

       1      6      6      1
o12 = R  <-- R  <-- R  <-- R  <-- 0
                                   
      0      1      2      3      4

o12 : ChainComplex

i13 : needsPackage "BettiCharacters";

i14 : A=action(RI2,G,Sub=>false)

o14 = ChainComplex with 6 actors

o14 : ActionOnComplex
\end{verbatim}

The action is defined with the option \texttt{Sub=>false}, which
allows passing group elements as square matrices rather than one-row
matrices of substitutions as in \S \ref{sec:example-1:-symmetric}.
Next we compute the character of the $G$-action on the resolution of
$I^2$ in homological degree 3.

\begin{verbatim}
i15 : character(A,3)

o15 = GradedCharacter{{24} => Character{1, 1, 1, 1, 1, 1}}

o15 : GradedCharacter
\end{verbatim}

As expected, we obtain a trivial character concentrated in degree 24.

The \texttt{BettiCharacters} package can also compute the characters
of a finite group on the graded components of a module.  Using the
package \texttt{SymbolicPowers} \cite{SymbolicPowersSource}, we can
directly establish that the character of $G$ on the graded component
of degree 21 in $I^{(2)} / I^2$ is trivial.

\begin{verbatim}
i16 : needsPackage "SymbolicPowers";

i17 : Is2 = symbolicPower(I,2);

o17 : Ideal of R

i18 : M = Is2 / I2;

i19 : B = action(M,G,Sub=>false)

o19 = Module with 6 actors

o19 : ActionOnGradedModule

i20 : character(B,21)

o20 = GradedCharacter{{21} => Character{1, 1, 1, 1, 1, 1}}

o20 : GradedCharacter
\end{verbatim}

\newcommand{\etalchar}[1]{$^{#1}$}
\def\cprime{$'$} \def\Dbar{\leavevmode\lower.6ex\hbox to 0pt{\hskip-.23ex
  \accent"16\hss}D} \def\Dbar{\leavevmode\lower.6ex\hbox to 0pt{\hskip-.23ex
  \accent"16\hss}D} \def\Dbar{\leavevmode\lower.6ex\hbox to 0pt{\hskip-.23ex
  \accent"16\hss}D} \def\Dbar{\leavevmode\lower.6ex\hbox to 0pt{\hskip-.23ex
  \accent"16\hss}D}

\end{document}